\documentclass{article}
\usepackage[francais,english]{babel}
\usepackage{amsmath,amsfonts,amsthm,amssymb,amscd}
\usepackage{color}

 \usepackage{verbatim}

\renewcommand{\Re}{\mathop{\rm Re}\nolimits}
\renewcommand{\Im}{\mathop{\rm Im}\nolimits}

\newcommand{\p}{\partial}

\newcommand{\C}{{\mathbb C}}
\newcommand{\R}{{\mathbb R}}

\newcommand{\la}{\lambda}

\newcommand{\ty}{\infty}
\newcommand{\ri}{\rightarrow}

\newcommand{\CC}{{\cal C}}

\newcommand{\UU}{{\cal U}}
\newcommand{\VV}{{\cal V}}

\newcommand{\lag}{\langle}
\newcommand{\rag}{\rangle}

\newcommand{\dd}{{\textup d}}

\theoremstyle{plain}
\newtheorem{theorem}{Theorem}[section]
\newtheorem{proposition}[theorem]{Proposition}

\newtheorem{condition}[theorem]{Condition}
\theoremstyle{remark}
\newtheorem{remark}[theorem]{Remark}

\newcommand{\de}{\delta}
\numberwithin{equation}{section}

\author{
Karine \textsc{Beauchard}\footnote{
CMLA, ENS Cachan, CNRS, UniverSud, 61, avenue du Pr�sident Wilson, F-94230 Cachan, FRANCE.
email: Karine.Beauchard@cmla.ens-cachan.fr}, 
Vahagn \textsc{Nersesyan}\footnote{
Laboratoire de Math\'ematiques de Versailles,
Batiment Fermat, 45, avenue des Etats-Unis,
F-78035 Versailles cedex, FRANCE.
email: 	nersesyan@math.uvsq.fr}
\thanks{The authors were partially supported by the ``Agence Nationale de la Recherche'' (ANR),
Projet Blanc C-QUID number BLAN-3-139579}
}

\title{Semi-global weak stabilization of bilinear Schr\"odinger equations}
\date{}

\begin{document}
\maketitle

\begin{abstract}
We consider a linear Schr\"odinger equation, on a bounded domain, with bilinear control,
representing a quantum particle in an electric field (the control). 
Recently, Nersesyan proposed explicit feedback laws and
proved the existence of a sequence of times $(t_n)_{n \in \mathbb{N}}$ for which
the values of the solution of the closed loop system
converge weakly in $H^2$ to the ground state.
Here, we prove the convergence of the whole solution, as $t \rightarrow + \infty$.
The proof relies on control Lyapunov functions and
an adaptation of the LaSalle invariance principle to PDEs.
\\

\textbf{R\'esum\'e:} Stabilisation faible semi-globale d'\'equations de Schr\"odinger bilin\'eaires.
\\
Nous consid\'erons une \'equation de Schr\"odinger lin\'eaire, sur un domaine born\'e,
avec un contr\^ole bilin\'eaire, mod\'elisant une particule quantique dans un champ \'electrique (la commande).
R\'ecemment, Nersesyan a propos\'e des lois de r\'etroaction explicites et
d\'emontr\'e l'existence d'une suite de temps $(t_n)_{n \in \mathbb{N}}$ auxquels
les valeurs de la solution du syst\`eme boucl\'e
convergent faiblement dans $H^2$ vers l'\'etat fondamental.
Ici, nous d\'emontrons la convergence de toute la solution, quand $t \rightarrow + \infty$.
La preuve repose sur des fonctions de Lyapunov et une adaptation du principe d'invariance de LaSalle 
aux EDP.

\end{abstract}

\newpage

\textbf{Version francaise abr\'eg\'ee: } On consid\`ere le syst\`eme
\begin{align} 
i\dot z   &= -\Delta z+V(x)z+ u(t) Q(x)z,\,\,\,\,x\in D,\label{E:hav1}\\
z\arrowvert_{\partial D}&=0,\label{E:ep1}\\
 z(0,x)&=z_0(x),\label{E:sp1}
\end{align}
o\`u $D\subset\R^m$ est un domaine born\'e    \`a bord lisse,
$V,Q\in C^\ty(\overline{D},\R)$ sont des fonctions donn\'ees,
$u$ est le contr\^ole et $z$ est l'\'etat. Il mod\'elise une particule quantique dans
un potentiel $V$ et un champ \'electrique $u$.

Notons $(e_{k,V})_{k \in \mathbb{N}^*}$ les vecteurs propres de l'op\'erateur $(-\Delta+V)$,
$(-\Delta+V) e_{k,V} = \lambda_{k,V} e_{k,V}$,
$P_{1,V}z:=z-\lag z,e_{1,V}\rag e_{1,V}$ la projection orthogonale de $L^2(D,\mathbb{C})$ sur 
$\text{Vect}\{ e_{k,V} , k\ge 2 \}$ et $S$ la sph\`ere unit\'e de $L^2(D,\mathbb{C})$.

Les lois de r\'etroaction explicites suivantes sont introduites dans \cite{N3}
\begin{equation}\label{E:fed3}
  {{ {u}}}(z):=-\delta\Im\Big[ \langle \alpha (-\Delta +V)   P_{1,V}(Qz ), (-\Delta
 +V)
 P_{1,V}z\rangle
 - \lag Qz,e_{1,V}\rag\lag e_{1,V},z\rag \Big],
 \end{equation} 
o\`u $\delta, \alpha>0$. Elles permettent de consid\'erer le syst\`eme boucl\'e
\begin{equation}
i\dot z    = -\Delta z+V(x)z+ u(z) Q(x)z,\,\,\,\,x\in D.\label{E:hav22}
\end{equation}
Rappelons la condition suivante, introduite \'egalement dans \cite{N3}.
\begin{condition}\label{C:pFR}
Les fonctions $V,Q\in C^\ty(\overline{D},\R) $ v\'erifient
\begin{enumerate}
\item [(i)] $\langle Qe_{1,V},e_{j,V}\rangle\neq0$
  pour tout $j\ge2$,
\item  [(ii)]$\la_{1,V}-\la_{j,V}\neq \la_{p,V}-\la_{q,V}$ pour tout  $j,p,q\ge1$ tels que
  $\{1,j\}\neq\{p,q\}$ et $j\neq 1$.
\end{enumerate}
\end{condition}

Dans \cite{N3}, Nersesyan d\'emontre que, sous la Condition \ref{C:pFR},
il existe une suite de temps $(t_n)_{n \in \mathbb{N}}$
auxquels la solution du syst\`eme boucl\'e converge faiblement dans $H^2$ vers l'\'etat fondamental:
$z(t_n) \rightharpoonup e_{1,V}$ dans $H^2$ quand $n \rightarrow + \infty$.
Dans cet article, nous d\'emontrons que toute la solution converge:
$z(t) \rightharpoonup e_{1,V}$ dans $H^2$ quand $t \rightarrow + \infty$.

\begin{theorem} \label{MainThmFr}
On suppose la Condition \ref{C:pFR} v\'erifi\'ee.
Soit $\UU_t$ la r\'esolvante du syst\`eme boucl\'e (\ref{E:hav22}), (\ref{E:ep1}).
Il existe un ensemble fini ou d\'enombrable $J\subset \R_{+}^*$ tel que, 
pour tout $z_0\in S\cap H^1_{0}\cap H^2$ n'appartenant pas � $\CC:=\{c e_{1,V}: c\in\C,|c|=1\}$,
il existe $\alpha^*:=\alpha^*(\|z_0\|_2)>0$ tel que, 
pour tout $\alpha \in (0,\alpha^*)-J$,
$\UU_{t}(z_0 )\rightharpoonup  \CC$, dans $H^2$, quand $t\ri\ty$.
 \end{theorem}  

La preuve du Th\'eor\`eme \ref{MainThmFr} repose sur le principe d'invariance de LaSalle et
se fait en deux \'etapes. Dans un premier temps, on v\'erifie que l'ensemble invariant 
coincide localement avec $\mathcal{C}$.
Dans un deuxi\`eme temps, on d\'emontre la convergence. Pour cela, on montre que
les seules valeurs d'adh\'erence possibles, pour la topologie faible $H^2$,
de la solution du syst\`eme boucl\'e, sont dans $\mathcal{C}$.
Consid\'erant une valeur d'adh\'erence faible $H^2$, 
$\mathcal{U}_{t_n}(z_0) \rightharpoonup z_\infty$,
on d\'emontre qu'elle appartient \`a $\mathcal{C}$,
en montrant qu'elle engendre une solution invariante, $u[ \mathcal{U}_t (z_\infty)] \equiv 0$.
Pour cela, on d\'emontre que $u[ \mathcal{U}_{t_n+t}(z_0)] \rightarrow 0$ quand $n \rightarrow + \infty$
pour presque tout $t \in [0,+\infty)$ et on justifie le passage \`a la limite $[n \rightarrow + \infty]$
dans le feedback.

\newpage

%OOOOOOOOOOOOOOOOOOOOOOOOOOOOOOOOOOOOOOOOOOOOOOOOOOOOOOOOOOOOOOOOOOOOOOOOOO
\section{Introduction}

%--------------------------------------------------------------------------
\subsection{The system}

We consider the system (\ref{E:hav1})-(\ref{E:sp1}) 
where $D\subset\R^m$ is a bounded domain with smooth boundary, 
$V,Q\in C^\ty(\overline{D},\R)$ are  given functions,  $u$ is the control, and $z$
is the state. It represents a quantum particle in a potential $V$, in an electric field $u$.
The following proposition establishes the well-posedness of system (\ref{E:hav1})-(\ref{E:sp1})
(see \cite{CW} for a proof).

\begin{proposition}\label{L:LD}
For any   $z_0\in H_0^1\cap H^2$ (resp. $z_0\in L^2$) and for any
  $u\in L^1_{loc}([0,\ty),\R)$
  problem
(\ref{E:hav1})-(\ref{E:sp1}) has a unique solution $z\in
C([0,\infty), H^1_0 \cap H^2)$ (resp. $z\in C([0,\infty),L^2)$).  
Furthermore, the resolving operator
$\UU_t(\cdot,u):L^2\rightarrow~L^2 $ taking $z_0$ to
    $z(t)$ satisfies the relation
\begin{align}
\|\UU_t(z_0,u)\|&=\|z_0\|,\,\,\, \forall t\ge0.\label{E:barev}
\end{align}
\end{proposition}

In all this article, 
$\|.\|$ (resp. $\|.\|_{s}$) denotes the usual norm on $L^2(D,\mathbb{C})$ 
(resp. $H^s(D,\mathbb{C})$, for every $s \in \mathbb{N}^*$).
$S$ is the $L^2(D,\mathbb{C})$-sphere and
$$\langle f , g \rangle := \int_D f(x) \overline{g(x)} dx, \forall f,g \in L^2(D,\mathbb{C}).$$

%-----------------------------------------------------------------------
\subsection{Bibliography}

We refer to \cite{KB-JMC, KB-CL, Chambrion-et-al, VaNe} for exact or approximate controllability
results for the system (\ref{E:hav1})-(\ref{E:sp1}), with open loop controls.
This article is concerned with closed loop controls:
we search explicit feedback laws, that asymptotically stabilize the ground state.

In \cite{MM-PR-GT}, the same question is addressed for ODE models.
The control design relies on control Lyapunov functions, and the
convergence proof relies on the LaSalle invariance principle.
This reference deals with the situation where 
the linearized system around the ground state
is controllable. The degenerate case is studied in \cite{BCMR}.

The goal of this article is to adapt the result of \cite{MM-PR-GT} to
PDE models. Indeed, the LaSalle invariance principle is a powerful tool
to prove the asymptotic stability of an equilibrium for a finite dimensional
dynamic system. However, using it for infinite dimensional systems is more difficult
(because closed and bounded subsets are not compact).

A first possible adaptation consists in proving \emph{approximate} convergence results,
as for example in \cite{KB-MM, MM}.
A second possible adaptation consists in proving a \emph{weak} convergence, as, for example, 
in \cite{BS} and in this article. 
A third possible adaptation consists in proving a strong convergence,
as for example in \cite{BAN-JMC}. In this case, one needs an
additional compactness property for the trajectories of the closed loop system.
Another strategy consists in designing strict Lyapunov functions,
as for example in \cite{JMC-BAN-GB}.

%----------------------------------------------------------------------
\subsection{Stabilization strategy}

Let us recall the stabilization strategy proposed in \cite{N3}.  
We introduce the Lyapunov function
\begin{equation*}
\VV(z): =\alpha \|( -\Delta+{V} )P_{1,V}z\| ^2+1-|\lag z,e_{1,V}
\rag |^2,\,\,\,\, z\in S\cap H^1_{0}\cap H^2,
\end{equation*}
where $\alpha>0$, 
$(e_{k,V})_{k \in \mathbb{N}^*}$ are the eigenvectors of the operator $-\Delta+V$,
$(-\Delta+V) e_{k,V} = \lambda_{k,V} e_{k,V}$
and   $P_{1,V}z:=z-\lag z,e_{1,V}\rag e_{1,V}$ is  the orthogonal projection in $L^2$ 
onto the closure of $\text{Span}\{ e_{k,V} , k\ge 2 \}$. 
Notice that $\VV(z)\ge0$ for all $ z\in S\cap H^1_{0}\cap H^2$ and 
$\VV(z)=0$ if and only if $z=ce_{1,V}, |c|=1$. 
For any $ z\in S\cap H^1_{0}\cap H^2$, we have
\begin{align*}
\VV(z)\ge\alpha \|( -\Delta+{V} )P_{1,V}z\| ^2\ge\frac{\alpha}{2}
\|  \Delta (P_{1,V}z)\| ^2-C_1\ge\frac{\alpha}{4} \|  \Delta z\|
^2-C_2,
\end{align*} where  $C_1$ and
 $C_2$ are positive constants.
Thus
\begin{equation}\label{E:sah}
C(1+\VV(z))\ge \|z\|_2
\end{equation}
for some constant $C>0$.
Following the ideas of \cite{BCMR, N3}, we wish to choose a feedback
law $u (\cdot)$ such that
\begin{equation*}
\frac{\dd}{\dd t}\VV(z(t))\le0
\end{equation*}for the  solution
$z(t)$ of (\ref{E:hav1})-(\ref{E:sp1}). Let us assume that $\Delta
z(t)\in H_0^1\cap H^2$ for all $t\ge0$. Using (\ref{E:hav1}), we get
 \begin{align*}
\frac{\dd}{\dd t}\VV(z(t))&=2\alpha\Re \Big[ \langle (-\Delta +V)
P_{1,V}{\dot{
z}}, (-\Delta +V)  P_{1,V}z\rangle \Big]
- 2 \Re \Big[ \lag\dot z,e_{1,V}\rag\lag e_{1,V},z\rag \Big] \nonumber\\
&= 2\alpha\Re \Big[ \langle  (-\Delta +V) P_{1,V}(i\Delta z -iVz-iu  Qz
), (-\Delta
+V)  P_{1,V}z\rangle \Big] \nonumber\\
&\quad-2\Re \Big[ \lag  i\Delta z -iVz-iuQz,e_{1,V}\rag\lag
e_{1,V},z\rag \Big].
\end{align*}
Integrating by parts and using the facts that  $V$ is real valued, $ P_{1,V}$ commutes with $-\Delta+V$    and
$$(-\Delta +V)  P_{1,V}z|_{\p D}=z|_{\p D}=e_{1,V}|_{\p D}=0,$$ we obtain
 \begin{align*}
 & 2\alpha\Re \Big[ \langle  
-i(-\Delta +V)^2  P_{1,V}z, (-\Delta +V)  P_{1,V}z\rangle \Big] 
-2\Re \Big[ \lag  i\Delta z -iVz ,e_{1,V}\rag\lag e_{1,V},z\rag \Big]\nonumber\\
&=  2\alpha\Re \Big[ \langle  - i\nabla (-\Delta +V)   P_{1,V}z, \nabla
(-\Delta +V) P_{1,V}z\rangle \Big] \nonumber\\
&\quad+2\alpha\Re \Big[ \langle
-iV(-\Delta +V) P_{1,V}z, (-\Delta +V)
P_{1,V}z\rangle \Big] \nonumber\\
&\quad+2\la_{1,V}\Re \Big[ \lag  i z
,e_{1,V}\rag\lag e_{1,V},z\rag \Big]=0.
\end{align*}
Thus
\begin{align*}
\frac{\dd}{\dd t}\VV(z(t)) = 2u\Im \Big[ \alpha\lag (-\Delta +V)
P_{1,V}(Qz), (-\Delta +V) P_{1,V}z\rangle -\lag Qz,e_{1,V}\rag\lag
e_{1,V},z\rag \Big] .
\end{align*}
Let us  take $u(z)$ defined by (\ref{E:fed3}) where $\de>0$. Then
\begin{equation}\label{E:fed35}
 \frac{\dd}{\dd t} \VV(z(t))=-\frac{2}{\delta}u^2(z(t)),
\end{equation}
thus $t \mapsto \VV(z(t))$ is not increasing and one may expect that
$z(t) \rightarrow \mathcal{C}:=\{ ce_{1,V}, c \in \mathbb{C} , |c|=1 \}$,
in some sense, when $t \rightarrow + \infty$. We consider the closed loop system (\ref{E:hav22}).
The following proposition ensures the well posedness of this system
and the validity of the computations performed above.
\begin{proposition} \label{WP}
For any $z_0\in H_0^1\cap H^2$   problem
 (\ref{E:hav22}), (\ref{E:ep1}), (\ref{E:sp1}) has a unique
 solution $z\in C([0,\ty),  H_0^1\cap H^2)$.
Moreover if $\Delta z_0\in H_0^1\cap H^2$, then, $\Delta z\in C([0,\ty),  H_0^1\cap H^2)$.
\end{proposition}

The local well-posedness and the regularity of the solution
of (\ref{E:hav22}), (\ref{E:ep1}), (\ref{E:sp1}) 
is standard (see \cite{CW}). From the
  construction of the feedback law $u$ it follows that a finite-time
  blow-up in $H^1_0 \cap H^2$ is impossible. Hence the solution is global in time.

%--------------------------------------------------------------------------
\subsection{Main result}

Let us introduce the following condition on the functions $V$ and $Q$.

\begin{condition}\label{C:p}
The functions $V,Q\in C^\ty(\overline{D},\R) $ are such that:
\begin{enumerate}
\item [(i)] $\langle Qe_{1,V},e_{j,V}\rangle\neq0$
  for all $j\ge2$,
\item  [(ii)]$\la_{1,V}-\la_{j,V}\neq \la_{p,V}-\la_{q,V}$ for all  $j,p,q\ge1$ such that
  $\{1,j\}\neq\{p,q\}$ and $j\neq 1$.
\end{enumerate}
\end{condition}

See the papers \cite{ PISI, VaNe,  PMMS} for the proof of genericity of this condition. 
The below theorem is the main result of this article

\begin{theorem}\label{T:stab}
Let $\UU_t$ be the resolving operator of the closed loop system (\ref{E:hav22}), (\ref{E:ep1}).
Under Condition \ref{C:p},  there is a finite  or countable set $J\subset \R_{+}^*$ such that for any
$\alpha\notin J$ and $z_0\in S\cap H^1_{0}\cap H^2$ with $0<\VV(z_0)<1$  we have 
  \begin{equation}\label{E:conv}
 \UU_{t}(z_0 )\rightharpoonup  \CC\,\,\,\,\text{in
$H^2$ as $t\ri\ty$},
 \end{equation} where $\CC:=\{c e_{1,V}: c\in\C,|c|=1\}.$
 \end{theorem}  

\begin{remark}
This Theorem proves the semi-global stabilization of the ground state.
Indeed, for every $z_0\in S\cap H^1_{0}\cap H^2$ such that $z_0 \notin \CC$,
one may chose $\alpha=\alpha(\|z_0\|_2)>0$ small enough so that the condition
$0<\mathcal{V}(z_0)<1$ is fulfilled.
\end{remark}

%OOOOOOOOOOOOOOOOOOOOOOOOOOOOOOOOOOOOOOOOOOOOOOOOOOOOOOOOOOOOOOOOOOOOOOOOOOOO
\section{Convergence proof}
\label{sec:Cv}

The first step of the proof consists in checking that the LaSalle invariance set
locally coincides with $\mathcal{C}$.
 
\begin{proposition} \label{Prop:LaS_Inv}
We assume Condition \ref{C:p}.
There exists a finite or countable set $J\subset \R_{+}^*$ such that, 
for every $\alpha\notin J$,
for every  $z_0\in {S\cap}H^1_{0}\cap H^2$ with $\lag z_0,e_{1,V} \rag\neq0 $ and $u(\UU_t(z_0))=0$ for all $t\ge0$, 
then  $z_0\in \CC$.  
\end{proposition}

This proposition is proved in \cite{N3}.
The second step of the proof consists in proving the convergence.
First, we need the following preliminary result.

\begin{proposition}\label{T:abc}
Let $\UU_t$ be the resolving operator of the closed loop system (\ref{E:hav22}), (\ref{E:ep1}). 
Let $z_n\in  H_0^1\cap H^2$ be such that 
$z_n\rightharpoonup z_\infty$ in $H^2$ and $z_n \rightarrow z_\infty$ in $H^1_0$.
For every $T>0$, there exists $N \subset (0,T)$ with zero Lebesgue measure such that
\begin{enumerate}
\item $\UU_{t}(z_n) \rightharpoonup \UU_t(z_\infty)$ in $H^2$ and
$\UU_{t}(z_n) \rightarrow \UU_t(z_\infty)$ in $H^1_0$, $\forall t \in (0,T)-N$,
\item $u[\UU_t(z_n)] \rightarrow u[\UU_t(z_\infty)]$, $\forall t \in (0,T)-N$.
\end{enumerate} 
\end{proposition}

\noindent \textbf{Proof of Proposition \ref{T:abc}:}
\\
\noindent \emph{First step: Let us show that, if
$z_n\in  H_0^1\cap H^2$, $z_n\rightharpoonup z_\infty$ in $H^2$
and $z_n \rightarrow z_\infty$ in $H^1_0$,
then $u(z_n)\rightarrow  u(z_\infty)$.}
Then, the second conclusion of Proposition \ref{T:abc} will be a consequence
of the first one.

Notice that (\ref{E:fed3}) and the fact that $Q$ is real valued imply  that
 \begin{align*}
u(z)=-\de\Im \Big[ \langle \alpha Q(-\Delta +V)     z  , (-\Delta
 +V)
  z\rangle \Big] +\tilde{u}(z)=\tilde{u}(z),
 \end{align*}where
\begin{equation} \label{u_reecrit}
\begin{array}{ll}
\tilde{u}(z)=&-\delta\Im \Big[ \langle \alpha (-\Delta +V)   P_{1,V}(Qz
), (-\Delta
 +V)
 (-\lag z,e_{1,V}\rag e_{1,V})\rangle
 \\&+ \langle \alpha (-\Delta +V)   (-\lag Qz,e_{1,V}\rag e_{1,V}), (-\Delta
 +V)
 z\rangle\\&+ \langle \alpha  (-\nabla Q\cdot \nabla z- z\Delta Q), (-\Delta
 +V)
 z\rangle\\&- \lag Qz,e_{1,V}\rag\lag e_{1,V},z\rag \Big].
\end{array}
\end{equation}
Thus, passing to the limit in the previous equality, we get $u(z_n)\rightarrow  u(z_\infty)$.
\\

\noindent \emph{Second step: Let us prove the first conclusion of Proposition \ref{T:abc}.}
Let $z_n\in  H_0^1\cap H^2$ be such that $z_n\rightharpoonup z_\infty$ in $H^2$ and
$z_n \rightarrow z_\infty$ in $H^1_0$. For $T>0$ define the Banach space   $W :=\{z\in C([0,T],H^1_0 \cap H^2)  \text{ such that } \dot z \in
 L^2([0,T], L^2 )\}$  endowed with the norm $\|z\|_W:=\|z\|_{ C([0,T],H^1_0 \cap H^2)}+\| \dot z\|_{L^2([0,T], L^2 )}$. 
The sequence of functions $(t \in [0,T] \mapsto \UU_t(z_n) )_{n \in \mathbb{N}}$ is bounded in
$W$, and the embedding $W \rightarrow L^2((0,T),H^1_0)$ is compact, by Theorem 5.1 in \cite{JLL}.
Let $Y \in L^2((0,T),H^1_0)$ and $\varphi$ an extraction such that
$$\UU_.(z_{\varphi(n)}) \rightarrow Y(.) \text{ in } L^2((0,T),H^1_0).$$
Thanks to the Lebesgue reciprocal theorem, one may assume that
\begin{equation} \label{CVH1pp}
\UU_t(z_{\varphi(n)}) \rightarrow Y(t) \text{ in } H^1_0, \forall t \in (0,T)-N,
\end{equation}
where $N \subset (0,T)$ has zero Lebesgue measure (otherwise take another extraction).

For every $t^* \in (0,T)-N$, the sequence $(\UU_{t^*}(z_{\varphi(n)}))_{n \in \mathbb{N}}$ is bounded in $H^2$
and its only possible weak $H^2$ limit is $Y(t^*)$ because of (\ref{CVH1pp}). 
Thus the whole sequence converges:
$\UU_{t^*}(z_{\varphi(n)}) \rightharpoonup Y(t^*)$ in $H^2$. Therefore, we have
\begin{equation} \label{CvH1H2pp}
\UU_t(z_{\varphi(n)}) \rightharpoonup Y(t) \text{ in } H^2 \text{ and }
\UU_t(z_{\varphi(n)}) \rightarrow Y(t) \text{ in } H^1_0, 
\forall t \in (0,T)-N.
\end{equation}
We deduce from the first step that
$$u[ \UU_t(z_{\varphi(n)}) ] \rightarrow u[Y(t)], \forall t \in (0,T)-N.$$

Let $A:=-\Delta+V$. We fix $t^* \in (0,T)-N$. For every $n$, we have
$$\UU_{t^*}(z_{\varphi(n)}) = e^{-iAt^*} z_{\varphi(n)}
+ i \int_0^{t^*} e^{-iA(t^*-s)} u[ \UU_s(z_{\varphi(n)})] Q \UU_s(z_{\varphi(n)}) ds.$$
Passing to the limit $[n \rightarrow + \infty]$ in $H^1_0$ in this equality, using the dominated convergence theorem and the continuity of $e^{-iAt}z$ with respect to  $z$ in $H^1_0$ norm, we get
$$Y(t^*)=e^{-iAt^*}z_\infty + i\int_0^{t^*} e^{-iA(t^*-s)} u[Y(s)] Q Y(s) ds.$$
Thus, $Y(t)=\UU_t(z_\infty), \forall t \in (0,T)-N$ (uniqueness of the solution of the closed loop system).

This proves that the sequence  $(t \in [0,T] \mapsto \UU_t(z_n) )_{n \in \mathbb{N}}$
has a unique adherence value in $L^2((0,T),H^1_0)$. Therefore the whole sequence converges
(i.e. one may take $\varphi = \text{Id}$). We deduce from (\ref{CvH1H2pp}) that 
$$\UU_t(z_n) \rightharpoonup \UU_t(z_\infty) \text{ in } H^2 \text{ and }
\UU_t(z_n) \rightarrow \UU_t(z_\infty) \text{ in } H^1_0, 
\forall t \in (0,T)-N.$$
\hfill $\Box$

\begin{remark}
The key point of this proof is that the feedback law $u(z)$ 
is well defined for $z$ strictly less regular that $H^2$
(see (\ref{u_reecrit}): formally, $z \in H^{3/2}$ is sufficient).
\end{remark}

\textbf{Proof of Theorem \ref{T:stab}:}
Let $J$ be as in Proposition \ref{Prop:LaS_Inv} and $\alpha \notin J$.
Let $z_0 \in S \cap H^1_0 \cap H^2$ with $0<\mathcal{V}(z_0)<1$.
Let us prove that the weak $H^2$ $\omega$-limit set of
$\{ \UU_t(z_0) ; t \geqslant 0 \}$ is contained in $\mathcal{C}$.

Let $z_\ty \in H^1_0 \cap H^2$ and $t_n\ri+\ty$ be such that 
$\UU_{t_n}(z_0) \rightharpoonup {z_\ty}$ in $H^2$. Let us show that ${ z_\ty} \in \CC$.
One may assume that $\UU_{t_n}(z_0) \rightarrow {  z_\ty}$ in $H^1_0$.

There exists an extraction $\varphi$ and a subset $N_1 \subset (0,+\infty)$
with zero Lebesgue measure such that
$$u[\UU_{t_{\varphi(n)}+t}] \rightarrow 0, \forall t \in (0,+\infty)-N_1.$$
Indeed, the sequence of functions
$(t \in (0,+\infty) \mapsto u[\UU_{t_n+t}(z_0)] )_{n \in \mathbb{N}}$
tends to zero in $L^2(0,+\infty)$ because
$t \mapsto u[\UU_t(z_0)]$ belongs to $L^2(0,+\infty)$ (see (\ref{E:fed35})).

Let $T \in [0,+\infty)$. Thanks to Proposition \ref{T:abc}, 
there exists $N \subset (0,T)$ with zero Lebesgue measure such that
$$u[\UU_{t_{\varphi(n)}+t}(z_0)] \rightarrow u[\UU_t(z_\infty)],
\forall t \in (0,T)-N.$$
The uniqueness of the limit ensures that
$u[\UU_t(z_\infty)]=0, \forall t \in (0,T)-[N \cup N_1]$.

Finally, the function $t \mapsto u[\UU_t(z_\infty)]$ is continuous
and vanishes on $(0,T)-[N\cup N_1]$, thus it vanishes on $[0,T]$.
This holds for every $T>0$, thus $u[\UU_t(z_\infty)]=0, \forall t \in [0,+\infty)$.
As $\VV({  z_\ty})\le\VV(z_0)<1,$ we have $\lag {  z_\ty},e_{1,V}\rag\neq0$. 
Thanks to Proposition \ref{Prop:LaS_Inv}, we get   ${ z_\ty} \in \mathcal{C}$.  $\hfill \Box$

%OOOOOOOOOOOOOOOOOOOOOOOOOOOOOOOOOOOOOOOOOOOOOOOOOOOOOOOOOOO
\section{Conclusion, open problems, perspectives}

We have proposed explicit feedback laws, that asymptotically stabilize the ground state,
for the system (\ref{E:hav1})-(\ref{E:sp1}).
To design the feedback laws, we have used control Lyapunov functions.
The convergence holds semi-globally in $H^2$ and for the weak $H^2$-topolopy.
The proof relies on a adaptation of the LaSalle invariance principle to PDEs.

Generalizations with different regularities are possible: 
with Lyapunov functions inspired by the $H^s$ distance to the target,
one may prove weak $H^s$ stabilization.
Generalization to other bilinear equations (for instance wave equations) are possible.

Our proof uses compact injections between Sobolev spaces on a bounded domain.
Thus, the stabilization question when such compact injections cannot be used is still an open problem.

Another open problem concerns the simultaneous stabilization of $N$ identical Schr\"odinger equations,
around $N$ different eigenstates, with only one closed loop control. Indeed, 
if we design feedback laws in the same way as in this
article, then, the LaSalle invariance set does not coincide with the target.
Thus, new ideas need to be introduced to tackle this problem.

The same question for non linear Schr\"odinger equations is also an open problem.

\end{document}